\documentclass[11pt,a4paper,twoside]{article}

\usepackage{amsmath,amssymb,amsthm}

\theoremstyle{plain}
\newtheorem{satz}{Theorem}[section]
\newtheorem*{namedtheorem}{\theoremname}
\newcommand{\theoremname}{testing}
\newenvironment{satzmitname}[1]{\renewcommand{\theoremname}{#1}
   \begin{namedtheorem}}
   {\end{namedtheorem}}
\newtheorem{lemma}[satz]{Lemma}
\newtheorem{prop}[satz]{Proposition}

\newenvironment{bew}{\begin{proof}}{\end{proof}}
\theoremstyle{definition}

\numberwithin{equation}{section}

\newcommand{\Q}{\ensuremath{\mathbb Q}}
\newcommand{\Z}{\ensuremath{\mathbb Z}}
\newcommand{\N}{\ensuremath{\mathbb N}}
\newcommand{\HG}{\ensuremath{\mathrm H}}
\newcommand{\HF}{\ensuremath{\widehat{\mathrm H}}}

\newcommand{\tp}{^\mathrm T}
\newcommand{\Gal}{\ensuremath{\operatorname{Gal}}}
\newcommand{\vcd}{\ensuremath{\operatorname{vcd}}}
\newcommand{\Hom}{\ensuremath{\operatorname{Hom}}}

\newcommand{\GL}{\ensuremath{\mathrm{GL}}}
\newcommand{\SL}{\ensuremath{\mathrm{SL}}}
\newcommand{\SP}{\ensuremath{\mathrm{Sp}}}

\newcommand{\quer}{\overline}

\newcommand{\de}{\mathfrak}

\newcommand{\seqmod}[3]{{#1\,\,\equiv\,\, #2\!\!\!\!\mod #3}}

\renewcommand{\geq}{\geqslant}
\renewcommand{\leq}{\leqslant}
\renewcommand{\epsilon}{\varepsilon}

\begin{document}
\title{On \hbox{$p$}-periodicity in the Farrell cohomology of
        \hbox{$\mathrm{Sp}(p-1,\mathbb{Z}[1/n])$}}
\author{Cornelia Minette Busch\thanks{Research partially supported by a
postdoctoral grant (SB2001--0138) from the Ministerio de Educaci\'on,
Cultura y Deporte of Spain.}\\
}
\date{November 2007}
\maketitle

\begin{abstract}
We describe isomorphism patterns in the $p$-primary part
of the Farrell cohomology ring $\HF^*(\SP(p-1,\Z[1/n]),\Z)$
for any odd prime $p$ and suitable integers $0\neq n\in\Z$,
where $\SP(p-1,\Z[1/n])$ denotes the group of symplectic
$(p-1)\/\times\/(p-1)$--matrices. Moreover, we determine
the precise $p$-period of this ring.
\end{abstract}

{\small{\hspace*{2mm}

\hspace*{2mm}
}}

%%%%%%%%%%%%%%%%%%%%%%%%%%%%%%%%%%%%%%%%%%%%%%%%%

\section{Introduction}

%%%%%%%%%%%%%%%%%%%%%%%%%%%%%%%%%%%%%%%%%%%%%%%%%

The Farrell cohomology is defined for any group $G$ with finite virtual
cohomological dimension ($\vcd G < \infty$). It is a generalization of
the Tate cohomology for finite groups. Let $\HF^*(G,\Z)$ denote the
Farrell cohomology of the group $G$ with coefficients in the ring $\Z$.
For a prime $p$ the $p$-primary part of $\HF^*(G,\Z)$ is written
$\HF^*(G,\Z)_{(p)}$. We then have
$$\HF^*(G,\Z)\cong\prod_p\HF^*(G,\Z)_{(p)}\,,$$
where $p$ ranges over the primes such that $G$ has $p$-torsion.
A group $G$ of finite virtual cohomological dimension has periodic
cohomology if for some $d\neq 0$ there is an element $u\in\HF^d(G,\Z)$
which is invertible in the ring $\HF^*(G,\Z)$. Cup product with $u$
then gives a periodicity isomorphism
$$\HF^i(G,M)\cong\HF^{i+d}(G,M)$$
for any $\Z G$-module $M$ and any $i\in\Z$. Similarly,
$G$ has $p$-periodic cohomology if for some $d\neq 0$
there is an element $u\in\HF^d(G,\Z)_{(p)}$ which is invertible
in the ring $\HF^*(G,\Z)_{(p)}$. We recall that if $G$ is a group with
$\vcd G <\infty$, then $G$ has $p$-periodic cohomology if and only if
every elementary abelian $p$-subgroup of $G$ has rank at most $1$.
For more details see Brown~\cite{brownb}.

Let $R$ be a commutative ring with $1$. The general linear group $\GL(n,R)$
is defined to be the multiplicative group of invertible
$n\/\times\/n$--matrices over $R$. The symplectic group $\SP(2n,R)$ over
the ring $R$ is the subgroup of matrices $Y\in\GL(2n,R)$ that satisfy
$$Y\tp J Y = J :=\begin{pmatrix}\phantom -0 & I \\ -I & 0 \end{pmatrix}\,,$$
where $I$ is the $n\/\times\/n$--identity matrix and $Y\tp$ denotes the
transpose of $Y$.
For any odd prime $p$ the symplectic group $\SP(p-1,\Z[1/n])$
over the ring $\Z[1/n]$, $0\neq n\in\N$, has finite virtual cohomological
dimension and $p$-periodic cohomology.

In section~\ref{isom} we determine isomorphisms between the
cohomology groups in the $p$-primary part of the cohomology ring
$\HF^i\bigl(\SP(p-1,\Z[1/n]),\Z\bigr)$.

\begin{satzmitname}{Theorem~\ref{hauptres isom}}
Let $p$ be an odd prime.
Let $n$ be such that $\Z[1/n][\xi]$ and $\Z[1/n][\xi +\xi^{-1}]$ are
principal ideal domains and moreover $p\mid n$.
Then for any $i\in\Z$
$$\HF^i\bigl(\SP(p-1,\Z[1/n]),\Z\bigr)_{(p)}\cong
  \HF^{i+b}\bigl(\SP(p-1,\Z[1/n]),\Z\bigr)_{(p)}$$
with $b=y$, the greatest odd divisor of $p-1$, if and only if for
each $j\mid y$ a prime $q\mid n$ exists with inertia
degree $f_q$ such that $j\mid \frac{p-1}{2f_q}$.
If no such $q$ exists, then $b=2y$.
\end{satzmitname}

The inertia degree $f_q$ of a prime $q\in\N$ is the multiplicative
order of $q$ in the field $\mathbb{F}_p$.
In section~\ref{period} we determine the periodicity isomorphisms in
the $p$-primary part of the cohomology.

\begin{satzmitname}{Theorem~\ref{hauptres periode}}
Let $n\in\Z$ be such that $\Z[1/n][\xi]$ and $\Z[1/n][\xi+\xi^{-1}]$ are
principal ideal domains and $p\mid n$. Then the $p$-period of the
Farrell cohomology ring
$$\HF^*\bigl(\SP(p-1,\Z[1/n]),\Z\bigr)$$
equals $2y$, where $y$ is the greatest odd divisor of $p-1$.
\end{satzmitname}

In fact the condition on the integer $n$ is not very restrictive
since $\Z[1/n][\xi]$ and $\Z[1/n][\xi+\xi^{-1}]$ are principal
ideal domains if and only if $q_1\cdot\ldots\cdot q_h$ divides $n$,
where $q_1,\dots,q_h\in\N$ are primes that depend on the prime $p$.
The integer $h$ is the class number of $\Q[\xi]$, i.e., the order of the ideal
class group of $\Z[\xi]$.
For primes $p$ with odd relative class number $h^-$ the $p$-period of
the Farrell cohomology ring $\HF^*(\SP(p-1,\Z),\Z)$ is $2y$, where $y$
is odd and $p-1=2^ry$ for some $r>0$, $r\in\Z$ (see \cite{BuschHFSp}).
The relative class number is $h^- :=h/h^+$,
where $h^+$ denotes the class number of $\Q[\xi +\xi^{-1}]$.

We use the following result of Brown (\cite{brownb}, Corollary X.7.4).
Let $G$ be a group with finite virtual cohomological dimension
such that each elementary abelian $p$-subgroup of $G$
has rank $\leq 1$. Then
\begin{equation}\label{KSBrown normalizers}
  \HF^*(G,\Z)_{(p)}\cong
  \prod_{P\in\mathfrak P}\HF^*\bigl(N(P),\Z\bigr)_{(p)}\,,
\end{equation}
where $\mathfrak P$ is a set of representatives of conjugacy classes
of subgroups $P$ of order $p$ in $G$ and $N(P)$ is the normalizer of $P$.
The symplectic group $\SP(p-1,\Z[1/n])$ that we are considering satisfies
this property.
Ash~\cite{ash} uses the isomorphism~\eqref{KSBrown normalizers} in order
to compute the Farrell cohomology of the group $\GL(n,\Z)$ with coefficients
in $\Z/p\Z$ for an odd prime $p$ and $p-1\leq n\leq 2p-3$. Naffah~\cite{nadim}
considers normalizers of subgroups of prime order in $\mathrm{PSL}(2,\Z[1/n])$
in order to compute the Farrell cohomology of $\mathrm{PSL}(2,\Z[1/n])$.
Glover and Mislin~\cite{GloMis} show corresponding results for the
outer automorphism group of the free group in the $p$-rank one case.
For the case $p=3$ see also the result of Adem and Naffah~\cite{adem+naffah}:
they consider the cohomology of the group $\SL(2,\Z[1/q])$, $q$ a prime.

It is well-known that if $N(P)/C(P)$ is a finite group
whose order is prime to $p$, then
\begin{equation}\label{cohom normalizer}
  \HF^*\bigl(N(P),\Z\bigr)_{(p)}\cong
  \bigl(\HF^*(C(P),\Z)_{(p)}\bigr)^{N(P)/C(P)}.
\end{equation}
In order to compute the $p$-period of $\HF^*(N(P),\Z)$,
we consider the action of $N(P)/C(P)$ on the centralizer
$C(P)$ and on $\HF^*(C(P),\Z)_{(p)}$. We already know the
structure of $N(P)/C(P)$ and
of $C(P)$ (see \cite{buschCCSpn}).

I would like to thank Carles Casacuberta for many valuable discussions.

%%%%%%%%%%%%%%%%%%%%%%%%%%%%%%%%%%%%%%%%%%%%%%%%%

\section{Subgroups of order $p$ in symplectic groups}

%%%%%%%%%%%%%%%%%%%%%%%%%%%%%%%%%%%%%%%%%%%%%%%%%

%-----------------------------------------------------
\subsection{Algebraic number theory}\label{Zahlentheorie}
%-----------------------------------------------------

The conjugacy classes of matrices of odd prime order $p$ in
$\SP(p-1,\Z[1/n])$ are related to some classes of ideals
in $\Z[1/n][\xi]$, where $\xi$ is a primitive $p$th root
of unity.

The ring $\Z[\xi]$ is the ring of integers of the cyclotomic field
$\Q(\xi)$ and $\Z[\xi +\xi^{-1}]$ is the ring of integers of the
maximal real subfield $\Q(\xi +\xi^{-1})$ of $\Q(\xi)$.
For any integer $0\neq n\in\Z$ we consider the ring $\Z[1/n]$
and the extensions $\Z[1/n][\xi]$ and $\Z[1/n][\xi +\xi^{-1}]$.
It is well-known that these are Dedekind rings. Let $q$
be a prime in $\Z$. The ideal $(q)\subset\Z[\xi+\xi^{-1}]$
can be written as a product of prime ideals
$\de q_1^+\cdots\de q_r^+$ in $\Z[\xi+\xi^{-1}]$. Consider
the prime ideals $\de q\subset\Z[\xi]$ over the prime ideal
$\de q^+\subset\Z[\xi+\xi^{-1}]$. The ideal $\de q^+$ satisfies
one of the following three properties.
\begin{enumerate}
\item The prime $\de q^+$ is inert: $\de q^+\Z[\xi]=\de q$ is
a prime ideal in $\Z[\xi]$ that lies over $q$.
\item The prime $\de q^+$ splits:
$\de q^+\Z[\xi]=\de q\quer{\de q}$, where
$\de q$ is a prime ideal in $\Z[\xi]$ that lies over $q$.
\item The ramified case: $\de p^+\Z[\xi]=\de p^2$, where
$\de p:=(1-\xi)$ is the only prime ideal in $\Z[\xi]$ that lies
over $p$. Moreover
$\de p^+:=((1-\xi)(1-\xi^{-1}))=\de p\quer{\de p}$
is the only prime ideal in $\Z[\xi+\xi^{-1}]$ that lies over $p$.
\end{enumerate}
The Galois group $G:=\Gal(\Q(\xi)/\Q)$,
resp.\  $G:=\Gal(\Q(\xi+\xi^{-1})/\Q)$,
acts transitively on the set of prime ideals $\de q$,
resp.\ $\de q^+$, that lie over the prime $q\in\Z$.
Some Galois automorphisms fix the prime ideals
$\de q\subset\Z[\xi]$ over $q$. These define the group
$$G_{\de q}=\{\gamma\in\Gal(\Q(\xi)/\Q)\mid \gamma(\de q)=\de q\}.$$
The order of $G_{\de q}$ is $f_q$, the inertia degree of $q$.
For more details see the book of Neukirch~\cite{neukirchschappa}.

%-----------------------------------------------------
\subsection{Centralizers and normalizers}
%-----------------------------------------------------

From now on $C(P)$ denotes the centralizer and $N(P)$ denotes
the normalizer of a subgroup $P$ of odd prime order $p$ in
$\SP(p-1,\Z[1/n])$.

%----------
\subsubsection{The centralizer}\label{C in Sp(p-1,Z[1/n])}
%----------

The centralizer of a subgroup of order $p$ in $\SP(p-1,\Z[1/n])$
is determined by the primes that lie over the primes that divide $n$.
Indeed we show in \cite{buschCCSpn}, Theorem~4.2, that if $n\in\Z$
is such that $\Z[1/n][\xi]$ and $\Z[1/n][\xi +\xi^{-1}]$ are
principal ideal domains, then the centralizer $C(P)$ of a subgroup
$P$ of order $p$ in $\SP(p-1,\Z[1/n])$ satisfies
$$C(P)\cong \Z/2p\Z\times\Z^{\sigma^+}\,.$$
Here $\sigma^+ =\sigma$ if $p\nmid n$, $\sigma^+=\sigma +1$ if
$p\mid n$ and $\sigma$ is the number of primes in $\Z[\xi +\xi^{-1}]$ that
split in $\Z[\xi]$ and lie over primes in $\Z$ that divide $n$.
This result is related to the fact that the centralizer $C(P)$ of $P$
is isomorphic to the kernel of the norm mapping
$$\begin{array}{rcl}
\Z[1/n][\xi]^* & \longrightarrow & \Z[1/n][\xi +\xi^{-1}]^* \\
x & \longmapsto & x\quer x.
\end{array}$$
%

%----------
\subsubsection{The quotient of the normalizer by the centralizer}\label{N/C in Sp(p-1,Z[1/n])}
%----------

Let $n\in\Z$ be such that $\Z[1/n][\xi]$ and $\Z[1/n][\xi +\xi^{-1}]$
are principal ideal domains and moreover $p\mid n$.
In \cite{buschCCSpn}, Theorem~4.1, we see that the
normalizer $N(P)$ and the centralizer $C(P)$ of a subgroup $P$
of order $p$ in $\SP(p-1,\Z[1/n])$ satisfy
$$N(P)/C(P)\cong\Z/j\Z\,,$$
where $j\mid p-1$, $j$ odd. Moreover, for each $j$ with $j\mid p-1$,
$j$ odd, a subgroup $P$ of order $p$ in $\SP(p-1,\Z[1/n])$
exists with $N(P)/C(P)\cong\Z/j\Z$.

%%%%%%%%%%%%%%%%%%%%%%%%%%%%%%%%%%%%%%%%%%%%%%%%%

\section{The Farrell cohomology}

%%%%%%%%%%%%%%%%%%%%%%%%%%%%%%%%%%%%%%%%%%%%%%%%%

Let $C(P)$ denote the centralizer and $N(P)$ the normalizer
of a subgroup $P$ of odd prime order $p$ in $\SP(p-1,\Z[1/n])$.
In this section we consider the $p$-primary part of the Farrell
cohomology ring $\HF^*\bigl(N(P),\Z\bigr)$. By \eqref{cohom normalizer}
we first determine the cohomology of the centralizer $C(P)$.
Then we describe the invariants under the action of the quotient
$N(P)/C(P)$ on $\HF^*\bigl(C(P),\Z\bigr)_{(p)}$.

%-----------------------------------------------------
\subsection{The Farrell cohomology of the centralizer}
%-----------------------------------------------------

\begin{prop}
Choose $n\in\Z$ such that $\Z[1/n][\xi]$ and $\Z[1/n][\xi +\xi^{-1}]$ are
principal ideal domains and $p\mid n$. Here $\xi$ is a primitive
$p$th root of unity.
Let $\sigma$ denote the number of primes in $\Z[\xi +\xi^{-1}]$
that split and lie over the primes in $\Z$ that divide $n$.
Then the Farrell cohomology ring of the centralizer
$C(P)$ of a subgroup $P$ of order $p$ in $\SP(p-1,\Z[1/n])$ is
\begin{align*}
  \HF^*\bigl(C(P),\Z\bigr) &\cong
  \HF^*\bigl((\Z/2p\Z)\times \Z^{\sigma +1},\Z\bigr) \\
  &\cong \Z/2p\Z[x,x^{-1}]\otimes \Lambda_\Z (e_0,\dots,e_{\sigma})\,,
\end{align*}
where $\deg (x)=2$ and $\deg (e_i)=1$, $i=0,\dots,\sigma$. In particular
$$\HF^i\bigl(C(P),\Z\bigr) \cong \bigl(\Z/2p\Z\bigr)^{2^{\sigma}}$$
and the $p$-primary part is
$$\HF^i\bigl(C(P),\Z\bigr)_{(p)} \cong \bigl(\Z/p\Z\bigr)^{2^{\sigma}}.$$
\end{prop}

\begin{bew}
It is a well-known result that
the entire ring of the Farrell cohomology of $(\Z/k\Z)\times \Z^l$,
$l\geq 1$, is
$$\HF^*\bigl((\Z/k\Z)\times \Z^l,\Z\bigr)\cong
  \Z/k\Z[x'',{x''}^{-1}]\otimes \Lambda_\Z (e_1,\dots,e_l)\,,$$
where $\deg (x'')=2$ and $\deg (e_i)=1$. In particular
$$\HF^i\bigl((\Z/k\Z)\times \Z^l,\Z\bigr)\cong
  \bigl(\Z/k\Z\bigr)^{2^{l-1}}.$$
A nice proof of this result is given in the thesis of
Naffah (\cite{nadim}, Proposition 5.10).
Now the assumption follows by \ref{C in Sp(p-1,Z[1/n])}.
For the $p$-primary part of the cohomology we have
\begin{align*}
  \HF^*\bigl(C(P),\Z\bigr)_{(p)} &\cong
  \bigl(\Z/2p\Z[x,x^{-1}]\otimes \Lambda_\Z (e_0,\dots,e_{\sigma})\bigr)_{(p)}\\
  &\cong \Z/p\Z[x',{x'}^{-1}]\otimes \Lambda_\Z (e_0,\dots,e_{\sigma})
\end{align*}
and herewith
$\HF^i\bigl(C(P),\Z\bigr)_{(p)} \cong \bigl(\Z/p\Z\bigr)^{2^{\sigma}}$.
\end{bew}

%----------
\subsubsection{The cup product}\label{cup product}
%----------

We have
$$\HF^*\bigl(C(P),\Z\bigr)_{(p)} \cong\Z/p\Z[x,x^{-1}]\otimes
\Lambda_{\Z/p\Z} (e_0,\dots,e_{\sigma})\,,$$
where $\deg (x)=2$ and $\deg (e_i)=1$, $i=0,\dots,\sigma$.
The cup product of
$x^k\otimes e, x^l\otimes e' \in \HF^*\bigl(C(P),\Z\bigr)_{(p)}$
is
\begin{align*}
(x^k\otimes e)\cdot (x^l\otimes e') &= (-1)^{\deg x^l \deg e} x^k x^l \otimes ee' \\
 &= x^{k+l}\otimes ee'
\end{align*}
because the degree of $x$ is even. In particular we get
$$(x^k\otimes e)\cdot (x^l\otimes e) = x^{k+l}\otimes ee = 0$$
if and only if $1\neq e\in\Lambda_{\Z/p\Z} (e_0,\dots,e_{\sigma})$
since in this case $ee=0$. Herewith the only invertible elements in
$\HF^*\bigl(C(P),\Z\bigr)_{(p)}$ are $x^k\otimes 1$, $k\in\Z$. Indeed
we have
$$(x^k\otimes 1)\cdot (x^{-k}\otimes 1)
= 1\otimes 1\in\HF^0\bigl(C(P),\Z\bigr)_{(p)}.$$
Moreover $\HF^*\bigl(C(P),\Z\bigr)_{(p)}$ is periodic of period
$2$ and the periodicity isomorphism is given by cup
product with $x\otimes 1$.

%-----------------------------------------------------
\subsection{An action on the Farrell cohomology of the centralizer}
%-----------------------------------------------------

Let $N(P)$ be the normalizer and $C(P)$ the centralizer of a subgroup
$P$ of odd prime order $p$ in $\SP(p-1,\Z[1/n])$. Choose $n\in\Z$ such that
$\Z[1/n][\xi]$ and $\Z[1/n][\xi +\xi^{-1}]$ are principal ideal domains
and moreover $p\mid n$. Here $\xi$ denotes a primitive $p$th root of unity.
In order to understand the action of $N(P)/C(P)$ on
$\HF^*(C(P),\Z)_{(p)}$ we recall how this quotient acts
on the centralizer $C(P)$. By \ref{C in Sp(p-1,Z[1/n])} the
sequence
$$C(P)\hookrightarrow \Z[1/n][\xi]^*\overset{N}{\rightarrow} \Z[1/n][\xi +\xi^{-1}]^*$$
is exact in $\Z[1/n][\xi]^*$. The norm $N$ is not surjective.
By \ref{N/C in Sp(p-1,Z[1/n])} the group
$N(P)/C(P)$ is isomorphic to a subgroup of the Galois group
$\Gal(\Q(\xi)/\Q)$:
$$N(P)/C(P) \hookrightarrow \Gal(\Q(\xi+\xi^{-1})/\Q) \hookrightarrow \Gal(\Q(\xi)/\Q).$$
The first embedding exists because the order of $N(P)/C(P)$ is odd.
Therefore the action of $N(P)/C(P)$ on the centralizer $C(P)$
is given by the action of $\Gal(\Q(\xi)/\Q)$ on the group of
units $\Z[1/n][\xi]^*$ and $N(P)/C(P)$ acts faithfully
on $C(P)$.

Now we determine the action of $N(P)/C(P)$ on $\HF^*\bigl(C(P),\Z\bigr)_{(p)}$.
We have
$$\HF^*\bigl(C(P),\Z\bigr)_{(p)} \cong
  \Z/p\Z[x,x^{-1}]\otimes \Lambda_\Z (e_0,\dots,e_{\sigma})\,,$$
where $x\in\HF^2(\Z/p\Z,\Z)$ and
$e_i\in\HG^1(\Z^{\sigma +1},\Z)=\Hom(\Z^{\sigma +1},\Z)$.

%----------
\subsubsection{The action on the first factor}\label{1.factor}
%----------

We know that $N(P)/C(P)$ is cyclic of order $j$, where $j\mid p-1$ and $j$
is odd. Since $N(P)/C(P)$ acts faithfully on $C(P)$, the action
of a generator of $N(P)/C(P)\cong\Z/j\Z$ on $\Z/p\Z[x,x^{-1}]$ is given
by $x\mapsto \mu x$, where $\mu\in(\Z/p\Z)^*$ is a primitive $j$th root of unity.
Then $x^j\mapsto \mu^j x^j = x^j$ and, in particular,
$$x^j \otimes 1\in \HF^{2j}\bigl(C(P),\Z\bigr)_{(p)}^{N(P)/C(P)} $$
is invertible in $\HF^{*}\bigl(N(P),\Z\bigr)_{(p)}$.

%----------
\subsubsection{The action on the second factor}\label{2.factor}
%----------

We have
\begin{align*}
\HF^*\bigl(C(P),\Z\bigr)_{(p)}
 &\cong
   \Z/p\Z[x,x^{-1}]\otimes \Lambda_\Z (e_0,\dots,e_{\sigma})\\
 &\cong
   \Z/p\Z[x,x^{-1}]\otimes_{\Z/p\Z} \Lambda_{\Z/p\Z} (e_0,\dots,e_{\sigma})
\end{align*}
and therefore we can consider the second factor to be the exterior
product of a $\Z/p\Z$--vector space.
The quotient $N(P)/C(P)\cong\Z/j\Z$ is isomorphic to a subgroup of the
Galois group $\Gal(\Q(\xi +\xi^{-1})/\Q)$ and the free abelian part of the
centralizer is given by the primes $\de p^+,\de q_1^+,\dots,\de q_\sigma^+$
in $\Z[\xi +\xi^{-1}]$ that are ramified or split in $\Z[\xi]$ and lie over the
primes that divide $n$. Any element $w\in C(P)$ can be written as
$$w=w'\epsilon_\de p^{m_0}\epsilon_1^{m_1}\dots\epsilon_\sigma^{m_\sigma}\,,$$
where $w'$ is the torsion part and $\epsilon_\de p\in\de p^+$,
$\epsilon_i\in\de q_i^+$, $i=1,\dots,\sigma$.
For a given basis we define a homomorphism
$$
\begin{array}{rcl}
\de p^+ \de q_1^+ \dots\de q_\sigma^+
& \longrightarrow
& \Z^{\sigma +1}
\\
\epsilon_\de p^{m_0}\epsilon_1^{m_1}\dots\epsilon_\sigma^{m_\sigma}
& \longmapsto
& (m_0,m_1,\dots,m_\sigma)
\end{array}
$$
and the dual elements $e_i\in\Hom(\Z^{\sigma +1},\Z)$, $i=0,\dots,\sigma$,
such that
$$e_i(m_0,\dots,m_\sigma)=\sum_{j=0}^\sigma \delta_{ij}m_j.$$
Let $E$ be the $\Z/p\Z$--vector space spanned by
$e_0,\dots,e_\sigma$.
The action of $\Gal(\Q(\xi +\xi^{-1})/\Q)$ on the primes
$\de p^+,\de q_1^+,\dots,\de q_\sigma^+$ defines an
action on $E$. The quotient $N(P)/C(P)$ acts as a subgroup
of $\Gal(\Q(\xi +\xi^{-1})/\Q)$. By the Herbrand unit theorem
the basis $\epsilon_{\de p},\epsilon_1,\dots,\epsilon_\sigma$
can be chosen such that the Galois group acts as a
permutation on this basis and herewith the group also
acts as a permutation on the basis $e_0,\dots,e_\sigma$
of $E$. Let $\gamma\in N(P)/C(P)\subseteq G=\Gal(\Q(\xi)/\Q)$
be a generator.
The group $\langle\gamma\rangle\subseteq G$
permutes the primes that lie over $q\mid n$. Therefore,
for each $q\mid n$ with $q$ prime and split, we have a
subspace $E_q\subseteq E:=\Z/p\Z(e_1,\dots,e_\sigma)$
that is invariant under the action of the Galois group.
Let $q\mid n$ be a prime such
that the prime $\de q^+\subset\Z[\xi +\xi^{-1}]$ that lies over $q$ splits,
i.e., $\de q^+\Z[\xi] = \de q\quer{\de q}$. The order of $\gamma$ is odd
and therefore $\gamma$ fixes $\de q$ (and $\quer{\de q}$)
if and only if $\gamma$ fixes $\de q^+$ and, moreover, the order $f_q$
of the group
$G_{\de q}=\{\widetilde\gamma\in G\mid \widetilde\gamma(\de q) = \de q\}$
satisfies $f_q\mid\frac{p-1}{2}$. The action of $N(P)/C(P)\cong\Z/j\Z$ on
$E_q$ is a permutation of order $c_{j,q}:=(\frac{p-1}{2f_q},j)$.
Therefore the eigenvalues of this action are $c_{j,q}$-th roots of
unity and $E_q$ is a direct sum of invariant subspaces $E_{j,q}$ of
dimension $c_{j,q}$. The characteristic polynomial of the
restriction of the action of a generator $\gamma\in N(P)/C(P)$
on the invariant subspaces $E_{j,q}$ is $1-x^{c_{j,q}}$ and the
charac\-teristic polynomial of the action of $\gamma$
on $E$ is a polynomial of the form
$$\prod_{q\mid n} (1 - x^{c_{j,q}})^{d_{j,q}}\,,$$
where $d_{j,q}$ is the number of subspaces of $E_q$ that are isomorphic
to $E_{j,q}$. By the definition of $c_{j,q}$ we have $c_{j,q}\mid j$
and herewith the eigenvalues of the action of $\gamma\in N(P)/C(P)$
on the space $E$ are $j$th roots of unity $\mu^k$, $0\leq k\leq j-1$.
The dimension of the invariant subspace $E_p$ is $1$ and $N(P)/C(P)$
acts trivially on $E_p$. Therefore $c_{j,p}=1$ and $d_{j,p}=1$.
For our purpose it is important which roots of unity occur as eigenvalues
but, as soon as it is nonzero, the multiplicity of those eigenvalues is
irrelevant.

%----------
\subsubsection{The action on the cohomology ring}\label{eigenvalue 1}
%----------

We consider the action of $N(P)/C(P)\cong \Z/j\Z$ on
$\HF^{2k+m}(C(P),\Z)_{(p)}$. We have seen in \ref{1.factor}
that the action of a generator $\gamma$ of $N(P)/C(P)$
on $x^k\in\Z/p\Z[x,x^{-1}]$ is given by multiplication with
$\mu^k$, where $\mu\in (\Z/p\Z)^*$ is a primitive $j$th root
of unity. If the same generator $\gamma$ acts on $e\in\Lambda^m E$
by multiplication with $\mu^l$, i.e., $e$ is an eigenvector to the
eigenvalue $\mu^l$, then $\gamma$ acts on
$$x^k\otimes e\in\HF^{2k+m}(C(P),\Z)_{(p)}$$
by multiplication with $\mu^l\mu^k$. This shows that
the element $x^k\otimes e$ is an eigenvector to the eigenvalue
$\mu^{l+k}$ of the action of $\gamma\in N(P)/C(P)$.
Since we are interested in the $N(P)/C(P)$--invariants of the $p$-primary
part of $\HF^i(C(P),\Z)$, we are searching for
the elements to the eigenvalue $1=\mu^{l+k}$, $l+k\in j\Z$.

%-----------------------------------------------------
\subsection{An example}
%-----------------------------------------------------

Let $p=7$, $\xi$ a primitive seventh root of unity and $n:=7q$,
where $q\in\Z$ is a prime such that the primes
$\de q_i^+\subset\Z[\xi +\xi^{-1}]$, $i=1,2,3$, that lie over $q$
split in $\Z[\xi]$, i.e., $f_q=1$. Since $f_q$ is the smallest
positive integer that satisfies $q^{f_q}\equiv 1\mod 7$, we
immediately see that $q:=29$ satisfies our condition and we
therefore choose $n:=7\cdot 29 = 203$.
Let $\de p^+ = ((1-\xi)(1-\xi^{-1}))$ be the prime in
$\Z[\xi+\xi^{-1}]$ over $p=7$.
The centralizer of a subgroup $P$ of order $7$ in
$\SP(6,\Z[1/203])$ is
$$C(P)\cong \Z/14\Z\times \Z^4\,.$$
We know that $N(P)/C(P)\cong\Z/j\Z$ with $j=1,3$. If $j=1$, then
$N(P)=C(P)$. We assume that we have chosen $P$ with
$N(P)/C(P)\cong\Z/3\Z$, i.e., $j=3$. Such a subgroup always
exists. A generator of the quotient $N(P)/C(P)$ acts as a
permutation: we choose the numbering of the $\de q_i$ such that
$\de q_1^+\mapsto\de q_2^+$,
$\de q_2^+\mapsto\de q_3^+$,
$\de q_3^+\mapsto\de q_1^+$.
We know that $\de p^+ \mapsto\de p^+$.
Then, by \ref{2.factor}, $\epsilon_0\in\de p$, $\epsilon_i\in\de q_i^+$,
and $e_0,e_i\in\Hom(\Z^4,\Z)$, $1\leq i \leq 3$, exist such that
the generator of $N(P)/C(P)$ acts
as a permutation on $E := \langle e_0, e_1, e_2, e_3\rangle$, i.e.,
$$\begin{array}{rclcrcl}
e_0 & \longmapsto & e_0\,, & \quad & e_2 & \longmapsto & e_3\,, \\
e_1 & \longmapsto & e_2\,, & \quad & e_3 & \longmapsto & e_1\,. \\
\end{array}$$
We see that $E=E_7\oplus E_{29}$, where $E_7=\langle e_0\rangle$
and $E_{29}=\langle e_1, e_2, e_3\rangle$ are the subspaces that are
invariant under the action of $N(P)/C(P)$.
We have $c_{3,29}=3$, $d_{3,29}=1$ and
$c_{3,7}= 1$, $d_{3,7}=1$.

Now we consider the Farrell cohomology of the centralizer:
$$\HF^*\bigl(C(P),\Z\bigr)_{(7)} \cong
  \Z/7\Z[x,x^{-1}]\otimes \Lambda_{\Z/7\Z} (e_0,e_1,e_2,e_3)\,,$$
where $x\in\HF^2(\Z/7\Z,\Z)$ and
$e_0,e_i\in\HG^1(\Z^4,\Z)=\Hom(\Z^4,\Z)$, $i=1\leq i\leq 3$. The
cohomology groups are $\Z/7\Z$--vector spaces:
\begin{equation}\label{HFi}
  \HF^i\bigl(C(P),\Z\bigr)_{(p)} \cong
  \sum_{\stackrel{0\leq m\leq 4}{\seqmod{m}{i}{2}}}
  \langle x^{\frac{i-m}{2}}\rangle\otimes \Lambda^m (e_0,\dots,e_3)\,,
\end{equation}
where $ \langle x^{\frac{i-m}{2}}\rangle$ is the $\Z/7\Z$--vector space
spanned by $x^{\frac{i-m}{2}}$.
The periodicity isomorhism is given by cup product with
the element $x\otimes 1\in \HF^2\bigl(C(P),\Z\bigr)_{(7)}$.

We consider the action of $N(P)/C(P)\cong \Z/3\Z$ on
$\HF^{2k+m}(C(P),\Z)_{(7)}$ since we are searching for the
invariants under this action. By \ref{eigenvalue 1} we first
determine the eigenspaces of $\Lambda^m (e_0,e_1,e_2,e_3)$,
$m=0,\dots,4$, under the action of a generator $\gamma$
of $\Z/3\Z$. The eigenvalues are the third roots of unity
$1,2,4\in(\Z/7\Z)^*$.
Then we choose the third root of unity $\mu = 2\in\Z/7\Z$ for the action of
$\Z/3\Z$ and get a $\Z/7\Z$--basis for any cohomology group:
$$
\begin{array}{l}
\HF^{0}(C(P),\Z)_{(7)}^{\Z/3\Z} = \\
\ \langle 1\otimes 1,\ x^{-1}\otimes (e_0e_1 + 4 e_0e_2 + 2 e_0e_3),\
x^{-1}\otimes (2 e_1e_2 + 3 e_1e_3 + e_2e_3) \rangle\,, \\[2mm]
\HF^{1}(C(P),\Z)_{(7)}^{\Z/3\Z} = \\
\ \langle 1\otimes e_0,\ 1\otimes (e_1+ e_2 + e_3),
\ x^{-1}\otimes (e_0e_1e_2 + 5 e_0e_1e_3 + 4 e_0e_2e_3)\rangle\,, \\[2mm]
\HF^{2}(C(P),\Z)_{(7)}^{\Z/3\Z} = \\
\ \langle 1\otimes (e_0e_1 + e_0e_2 + e_0e_3),
\ 1\otimes (e_1e_2 - e_1e_3 + e_2e_3)\rangle\,, \\[2mm]
\HF^{3}(C(P),\Z)_{(7)}^{\Z/3\Z} = \\
\ \langle x\otimes (e_1 + 2 e_2 + 4 e_3),
\ 1\otimes (e_1e_2e_3),
\ 1\otimes (e_0e_1e_2 - e_0e_1e_3 + e_0e_2e_3)\rangle\,, \\[2mm]
\HF^{4}(C(P),\Z)_{(7)}^{\Z/3\Z} = \\
\ \langle x\otimes (e_0e_1 + 2 e_0e_2 + 4 e_0e_3),
\ x\otimes (e_1e_2 + 3 e_1e_2 + 2e_2e_3),
\ 1\otimes e_0e_1e_2e_3 \rangle\,, \\[2mm]
\HF^{5}(C(P),\Z)_{(7)}^{\Z/3\Z} = \\
\ \langle x^2\otimes (e_1 + 4 e_2 + 2 e_3),
\ x\otimes (e_0e_1e_2 + 3 e_0e_1e_3 + 2 e_0e_2e_3)\rangle\,.
\end{array}
$$
If we choose $\mu = 4\in\Z/7\Z$, we get other generators for the
cohomology rings, but the cohomology groups are isomorphic.
Cup product with the invertible element
$$x^3\otimes 1\in\HF^{6}(C(P),\Z)_{(7)}^{\Z/3\Z}$$
yields a periodicity isomorphism of degree $6$.
But we see that more cohomology groups are
isomorphic as $\Z/7\Z$--vector spaces. Indeed
$$\HF^{i}(C(P),\Z)_{(7)}^{\Z/3\Z}\cong \HF^{i+3}(C(P),\Z)_{(7)}^{\Z/3\Z}$$
for any $i\in\Z$. In this article we determine under which
conditions this isomorphism exists. We first explain the general
discussion on this example.

The dimension of the eigenspace of $\Lambda^m E$,
$m=0,\dots,4$, to the eigenvalue $\mu^l$, $l=0,1,2$,
$\mu\in (\Z/7\Z)^*$, is given by the coefficients $D_m[l]$
of $t^m \mu^l$ in the polynomial
\begin{align*}
L(t,\mu) &:= \sum_{m,l}D_m[l]t^m \mu^l =
  (1+t)(1+t\mu)(1+t\mu^2)(1+t) \\
  &= 1 + 2 t + t \mu + t \mu^2 + 2 t^2 + 2 t^2 \mu + 2 t^2 \mu^2
     + 2 t^3 + t^3 \mu + t^3 \mu^2 + t^4 .
\end{align*}
The variable $t$ counts the degree of the elements.
By the isomorphism~\eqref{HFi} we get
$$\dim \Bigl(\HF^i\bigl(C(P),\Z\bigr)_{(p)}^{N(P)/C(P)}\Bigr) =
  \sum_{\stackrel{0\leq m\leq 4}{\seqmod{m}{i}{2}}}
  D_m\Bigl[\frac{m-i}{2}\Bigr]$$
and this formula yields the dimensions of the cohomology groups
that we explicitly determined before.
In the next section we make the general discussion
of the arguments that we presented here.

%%%%%%%%%%%%%%%%%%%%%%%%%%%%%%%%%%%%%%%%%%%%%%%%%

\section{Isomorphisms in the cohomology ring}\label{isom}

%%%%%%%%%%%%%%%%%%%%%%%%%%%%%%%%%%%%%%%%%%%%%%%%%

Let $N(P)$ denote the normalizer and $C(P)$ the centralizer of a
subgroup $P$ of order $p$ of $\SP(p-1,\Z[1/n])$.
Let $n\in\Z$ be such that $\Z[1/n][\xi]$ and $\Z[1/n][\xi +\xi^{-1}]$ are
principal ideal domains and moreover $p\mid n$.

\begin{prop}\label{Iso von HF(N,Z)}
Let $P$ be such that $N(P)/C(P)\cong \Z/j\Z$ for a fixed odd $j>0$
with $j\mid p-1$. Then for any $i\in\Z$
$$\HF^i\bigl(N(P),\Z\bigr)_{(p)}\cong \HF^{i+b_j}\bigl(N(P),\Z\bigr)_{(p)}$$
with $b_j=j$ if and only if a prime $q\mid n$ exists with inertia
degree $f_q$ such that $j\mid\frac{p-1}{2f_q}$.
If no such $q$ exists, then $b_j=2j$.
\end{prop}

\begin{bew}
We consider the action of $N(P)/C(P)\cong \Z/j\Z$ on
the $\Z/p\Z$--vector space $\HF^{2k+m}(C(P),\Z)_{(p)}$.
By \ref{eigenvalue 1} we are searching for elements
$$x^k\otimes e\in\HF^{2k+m}\bigl(C(P),\Z\bigr)_{(p)}$$
to the eigenvalue $1$. If a generator
$\gamma\in N(P)/C(P)$ acts on
$x\otimes 1\in\HF^{2k+m}(C(P),\Z)_{(p)}$ by
multiplication with $\mu$ and on $e\in\Lambda^m E$ by
multiplication with $\mu^l$, i.e., $e$ is an eigenvector
to the eigenvalue $\mu^l$, then $\gamma$ acts on
$x^k\otimes e\in\HF^{2k+m}(C(P),\Z)_{(p)}$ by
multiplication with $\mu^k\mu^l$. We first
consider the eigenspaces of $\Lambda^m E$ under
the action of $N(P)/C(P)\cong \Z/j\Z$. Since this
quotient acts as a permutation on the space $E$
spanned by $\{e_0,e_1,\dots,e_\sigma\}$, the
eigenvalues are $j$th roots of unity $\mu^l\in(\Z/p\Z)^*$,
$l=0,\dots ,j-1$, and a basis $\{e'_0, e'_1,\dots, e'_\sigma\}$
of eigenvectors exists for $E$.
The elements of a basis of eigenvectors of the space
$\Lambda^m E$ are the products $e'_{i_1}\cdots e'_{i_m}$,
$i_1<\ldots < i_m$, where the eigenvalue of the product
equals the product of the eigenvalues.
Therefore the dimension of the eigenspace of $\Lambda^m E$
to the eigenvalue $\mu^l$ is given by the coefficients $D_m[l]$ of
$t^m\mu^l$ in the polynomial
$$L(t,\mu):= \sum_{m,l}D_m[l]t^m \mu^l =
  \prod_{\stackrel{q\mid n\text{ splits}}{\text{or }q=p\mid n}}
    \Bigl(\ %
      \prod_{k=1}^{c_{j,q}}\bigl(1+t\mu^{k\frac{j}{c_{j,q}}}\bigr)
    \ \Bigr)^{d_{j,q}}.$$
The variable $t$ in $L(t,\mu)$ counts the degree of the elements.
We have
$$\HF^i\bigl(C(P),\Z\bigr)_{(p)} =
  \sum_{\stackrel{0\leq m\leq \sigma +1}{\seqmod{m}{i}{2}}}
  \langle x^{\frac{i-m}{2}}\rangle\otimes \Lambda^m (e_0,\dots,e_\sigma)\,,$$
where $ \langle x^{\frac{i-m}{2}}\rangle$ is the $\Z/p\Z$--vector space
spanned by $x^{\frac{i-m}{2}}$.
We get
$$\dim \Bigl(\HF^i\bigl(C(P),\Z\bigr)_{(p)}^{N(P)/C(P)}\Bigr) =
  \sum_{\stackrel{0\leq m\leq \sigma +1}{\seqmod{m}{i}{2}}}
  D_m\Bigl[\frac{m-i}{2}\Bigr].$$
Consider the polynomial
$$L(z,z^{-2}) = \sum_l a_l z^l =
  \prod_{\stackrel{q\mid n\text{ splits}}{\text{or }q=p\mid n}}
    \Bigl(
      \prod_{k=1}^{c_{j,q}}\bigl(1+z^{1-2k\frac{j}{c_{j,q}}}\bigr)
    \Bigr)^{d_{j,q}}.
$$
Herewith we get
$$\dim \Bigl(\HF^i\bigl(C(P),\Z\bigr)_{(p)}^{N(P)/C(P)}\Bigr) =
   \sum_{\seqmod{l}{i}{2j}} a_l.$$
If $q\mid n$ exists with $c_{j,q}=j$, then the product
$$\prod_{k=1}^j 1+z^{1-2k}$$
is a factor of the polynomial $L(z,z^{-2})$ and for
$k=\frac{j+1}{2}$ ($j$ is odd) we get the factor $1+z^{-j}$
of $L(z,z^{-2})$. By Lemma~\ref{Polynomrechnung}
$$\dim \Bigl(\HF^i\bigl(C(P),\Z\bigr)_{(p)}^{N(P)/C(P)}\Bigr)
  = \dim \Bigl(\HF^{i+j}\bigl(C(P),\Z\bigr)_{(p)}^{N(P)/C(P)}\Bigr)$$
if and only if $1+ z^{j+2kj}$ is a factor of the polynomial
$L(z,z^{-2})$ for some $k\in\Z$. This happens if and only if $q$ exists
with $c_{j,q}=j$, i.e., if and only if $q$ exists with inertia degree
$f_q$ such that $j\mid \frac{p-1}{2f_q}$. The cohomology
groups that have the same dimension are isomorphic as
$\Z/p\Z$--vector spaces, but the isomorphism is not
always a periodicity isomorphism.
\end{bew}

\begin{lemma}\label{Polynomrechnung}
Let $\sum_{l\in \Z}a_l z^l$ be a polynomial with coefficients
in $\Z$. Then
$$\sum_{\seqmod{l}{i}{2j}}a_l = \sum_{\seqmod{l}{i+j}{2j}}a_l$$
if and only if $1+z^{-j}$ (or $1+z^{j+2kj}$ for some $k\in\Z$)
is a factor of the polynomial $\sum_{l\in \Z}a_l z^l$.
\end{lemma}

\begin{bew}
Consider the polynomial $f(z)=\sum_{l\in \Z}a_l z^l$, $a_l\in\Z$.
If $g(z)=\sum_{l\in \Z}a'_l z^l$, $a'_l\in\Z$, exists with
$f(z)=(1+z^{\pm j})g(z)$, then
$$f(z) = g(z) + z^{\pm j} g(z)
  = \sum_{l\in\Z} a'_l z^l + \sum_{l\in\Z} a'_l z^{l\pm j}$$
and herewith
\begin{align*}
\sum_{\seqmod{l}{i}{2j}} a_l
 &= \sum_{\seqmod{l}{i}{2j}} a'_l + \sum_{\seqmod{l}{i+j}{2j}} a'_l
  = \sum_{\seqmod{l}{i}{j}} a'_l \\
 &= \sum_{\seqmod{l}{i+j}{2j}} a_l\,.
\end{align*}
For the other direction we first consider the special case
$$\sum_{\seqmod{l}{0}{2j}}a_l = \sum_{\seqmod{l}{j}{2j}}a_l\,.$$
The value of the polynomial
$$\sum_{\seqmod{l}{0}{2j}}a_l (z^j)^{l/j} + \sum_{\seqmod{l}{j}{2j}}a_l (z^j)^{l/j}$$
is $0$ in $z^j=-1$. Therefore $(1+z^j)$ divides the polynomial.
The cases
$$\sum_{\seqmod{l}{i}{2j}}a_l = \sum_{\seqmod{l}{i+j}{2j}}a_l$$
are analogous.
The assumption now follows by an addition.
\end{bew}

\begin{satz}\label{hauptres isom}
Let $p$ be an odd prime.
Let $n$ be such that $\Z[1/n][\xi]$ and $\Z[1/n][\xi +\xi^{-1}]$ are
principal ideal domains and moreover $p\mid n$.
Then
$$\HF^i\bigl(\SP(p-1,\Z[1/n]),\Z\bigr)_{(p)}\cong
  \HF^{i+b}\bigl(\SP(p-1,\Z[1/n]),\Z\bigr)_{(p)}$$
for any $i\in\Z$, with $b=y$, the greatest odd divisor of $p-1$,
if and only if for each $j\mid y$ a prime $q\mid n$ exists with
inertia degree $f_q$ such that $j\mid \frac{p-1}{2f_q}$.
If no such $q$ exists, then $b=2y$.
\end{satz}

\begin{bew}
If $P$ is a subgroup of order $p$ in $\SP(p-1,\Z[1/n])$
that satisfies $N(P)/C(P)\cong \Z/j\Z$,
then we know by the proof of Proposition~\ref{Iso von HF(N,Z)} that
$$\HF^i\bigl(N(P),\Z\bigr)_{(p)}\cong \HF^{i+b_j}\bigl(N(P),\Z\bigr)_{(p)}$$
with $b_j=j$ if and only if a prime $q\mid n$ exists with inertia
degree $f_q$ such that $j\mid\frac{p-1}{2f_q}$.
If no such $q$ exists, then $b_j=2j$.
In order to determine the degree $b$ in our assumption, we let $P$ run
through the sets of conjugacy
classes of subgroups of order $p$ in $\SP(p-1,\Z[1/n])$.
Then the order $j$ of $N(P)/C(P)$ runs through the odd divisors of $y$.
By the isomorphism~\eqref{KSBrown normalizers}, the degree
$b$ is the least common multiple of the $b_j$ and this is $b=y$
if all the $b_j$ are odd and $b=2y$ if one of these numbers
is even. This proves the theorem.
\end{bew}

%%%%%%%%%%%%%%%%%%%%%%%%%%%%%%%%%%%%%%%%%%%%%%%%%

\section{The $p$-periodicity}\label{period}

%%%%%%%%%%%%%%%%%%%%%%%%%%%%%%%%%%%%%%%%%%%%%%%%%

Let $N(P)$ denote the normalizer and $C(P)$ the centralizer of a
subgroup $P$ of order $p$ of $\SP(p-1,\Z[1/n])$.
Let $n\in\Z$ be such that $\Z[1/n][\xi]$ and $\Z[1/n][\xi +\xi^{-1}]$ are
principal ideal domains and moreover $p\mid n$.

%----------
\subsection{The $p$-period of $\HF^*\bigl(N(P),\Z\bigr)$}
%----------

\begin{prop}\label{p-Periode von HF(N,Z)}
Let $P$ be such that $N(P)/C(P)\cong \Z/j\Z$ for a fixed odd $j>0$
with $j\mid p-1$. Then the periodicity isomorphism in
$$\HF^*\bigl(N(P),\Z\bigr)_{(p)}$$
is given by cup product
with $x^j\otimes 1\in\HF^{2j}\bigl(N(P),\Z\bigr)_{(p)}$ and
the period is $2j$.
\end{prop}

\begin{bew}
By \ref{cup product} the element
$$x\otimes 1\in\HF^2\bigl(C(P),\Z\bigr)_{(p)}$$
is invertible in the cohomology ring
and cup product with $x\otimes 1$ yields the periodicity
isomorphism. We know by \ref{1.factor} that the action of
a generator of $N(P)/C(P)$ on $x\otimes 1$ is given by
multiplication with a primitive $j$th root of unity
$\mu\in(\Z/p\Z)^*$. Therefore, by \ref{cup product},
$$x^j\otimes 1
\in\HF^{2j}\bigl(C(P),\Z\bigr)_{(p)}^{N(P)/C(P)}
\cong\HF^{2j}\bigl(N(P),\Z\bigr)_{(p)}\,$$
is invertible and cup product with $x^j\otimes 1$ yields the
periodicity isomorphism. The period is $2j$.
\end{bew}

%----------
\subsection{The $p$-period of $\HF^*(\SP(p-1,\Z[1/n]),\Z)$}
%----------

\begin{satz}\label{hauptres periode}
Let $n\in\Z$ be such that $\Z[1/n][\xi]$ and $\Z[1/n][\xi+\xi^{-1}]$ are
principal ideal domains and $p\mid n$. Then the $p$-period of the
Farrell cohomology ring
$$\HF^*\bigl(\SP(p-1,\Z[1/n]),\Z\bigr)$$
equals $2y$, where $y$ is the greatest odd divisor of $p-1$.
\end{satz}

\begin{bew}
If $P$ is a subgroup of order $p$ in $\SP(p-1,\Z[1/n])$
that satisfies $N(P)/C(P)\cong \Z/j\Z$,
then we know by the proof of Proposition~\ref{p-Periode von HF(N,Z)}
that the periodicity isomorphism of the corresponding factor
$$\HF^*\bigl(N(P),\Z\bigr)_{(p)}$$
in \eqref{KSBrown normalizers} is given by $x^j\otimes 1$ and
therefore the period of this factor equals $2j$.
If $P$ runs through the sets of conjugacy classes of subgroups of order $p$
in $\SP(p-1,\Z[1/n])$, then the order $j$ of $N(P)/C(P)$ runs through
the odd divisors of $y$. Therefore the least common multiple of the
$j$ is $y$. If $x^j\otimes 1$ is invertible in
$\HF^*\bigl(N(P),\Z\bigr)_{(p)}$, then $x^y\otimes 1$ is also invertible,
because $y$ is a multiple of $j$.
Now, by \eqref{KSBrown normalizers}, the
periodicity isomorphism of
$$\HF^*\bigl(\SP(p-1,\Z[1/n]),\Z\bigr)_{(p)}$$
is given by cup product with $x^y\otimes 1$ and therefore the
$p$-period is $2y$.
\end{bew}

%%%%%%%%%%%%%%%%%%%%%%%%%%%%%%%%%%%%%%%%%%%%%%%%%

\vspace*{1cm}

\begin{minipage}[t]{12cm}
Cornelia Minette Busch

Katholische Universit\"at Eichst\"att--Ingolstadt

MGF

D-85071 Eichst\"att

Germany
\end{minipage}


\begin{thebibliography}{99}
\bibitem{adem+naffah}
        A.~Adem, N.~Naffah,
        \emph{On the cohomology of ${SL}_2(\mathbb{Z}[1/p])$},
        {G}eometry and cohomology in group theory,
        London {M}ath.\ {S}oc.\ {L}ect.\ {N}ote Ser.\ 252, 1--9,
        Cambridge {U}niversity {P}ress (1998).

\bibitem{ash}
        A.~Ash,
        \emph{Farrell cohomology of ${GL}(n,\mathbb{Z})$},
        Israel J. Math. 67 (1989), 327--336.

\bibitem{brownb}
        K.~S.~Brown,
        \emph{Cohomology of groups}, Graduate Texts in Mathematics 87, Springer 1982.

\bibitem{BuschHFSp}
        C.~M.~Busch,
        \emph{The Farrell cohomology of \hbox{$\mathrm{Sp}(p-1,\mathbb{Z})$}},
        Doc.\ Math.\ 7 (2002), 239--254.

\bibitem{buschCCSpn}
        C.~M.~Busch,
        \emph{Conjugacy classes of \hbox{$p$}-torsion in symplectic
        groups over ${S}$-integers}, New York J. Math. 12 (2006), 169--182.

\bibitem{GloMis}
        H.~H.~Glover, G.~Mislin,
        \emph{On the $p$-primary cohomology of \hbox{$Out(F_n)$} in the $p$-rank
        one case},
        J. Pure Appl.\ Algebra 153 (2000), 45--63.

\bibitem{nadim}
        N.~Naffah,
        \emph{On the integral {F}arrell cohomology ring of ${PSL}_2(\mathbb{Z}[1/n])$},
        Diss. {E}{T}{H} {N}o. 11675, ETH Z\"urich, 1996.

\bibitem{neukirchschappa}
        J.~Neukirch,
        \emph{Algebraic number theory},
        Grundlehren der mathematischen Wissenschaften 322, Springer 1999.
\end{thebibliography}
\end{document}